\documentclass[reqno,11pt]{amsart}
\DeclareMathOperator*{\esssup}{\mathrm{ess\,sup}}

\usepackage{amssymb,amsthm}
\usepackage{amsmath}
\usepackage{amsfonts}
\usepackage{dsfont}

\usepackage[a4paper,  margin=3cm]{geometry}

\usepackage[active]{srcltx}
\makeatletter\@addtoreset{equation}{section}\makeatother

\newtheorem{theorem}{Theorem}[section]

\newtheorem{lemma}[theorem]{Lemma}

\newtheorem{assumption}[theorem]{Assumption}

\newtheorem{remark}[theorem]{Remark}
\numberwithin{equation}{section}

\title{Self-regulation in the Bolker-Pacala model}

\author{Yuri  Kondratiev}
\address{Fakut\"at f\"ur Mathematik, Universit\"at Bielefeld, Bielefeld D-33615, Germany and Interdisciplinary Center
for Complex Systems, Dragomanov University, Kyiv, Ukraine}
\email{kondrat@math.uni-bielefeld.de}

\author{ Yuri  Kozitsky}

\address{Instytut Matematyki, Uniwersytet Marii Curie-Sk{\l}odowskiej, 20-031 Lublin, Poland}
\email{jkozi@hektor.umcs.lublin.pl}


\begin{document}


\begin{abstract}

The Markov dynamics is studied of an infinite system of point
entities placed in $\mathds{R}^d$, in which the constituents
disperse and die, also due to competition. Assuming that the
dispersal and competition kernels are continuous and integrable we
show that the evolution of states of this model preserves their
sub-Poissonicity, and hence the local self-regulation (suppression
of clustering) takes place. Upper bounds for the correlation
functions of all orders are also obtained for both long and short
dispersals, and for all values of the intrinsic mortality rate.

\end{abstract}

\maketitle

\section{Introduction}


An actual task of applied mathematics is the description of the
dynamics of large systems of living entities, see \cite{BP1,Mu,O}.
This relates, in particular, to the systems in which the dynamics
amounts to the appearance (birth) and disappearance (death) of the
constituents. The disappearance of a given entity caused by its
interaction with the rest of the community is interpreted as the
result of  \emph{competition}.

In the simplest birth-and-death models, the state space is
$\mathds{N}_0 :=\mathds{N}\cup \{0\}$. Then the only observed result
of the trade-off between the appearance and disappearance is the
dynamics of the number of entities in the population. The theory of
such models goes back to A. Kolmogorov and W. Feller, see
\cite[Chapter XVII]{Feller}. In this theory, the time evolution of
the probability of having $n$ entities is obtained from the
Kolmogorov equation with a tridiagonal infinite matrix containing
birth and death rates. An important generalization here is to place
the entities in a continuous habitat, $\mathcal{H}$, usually a
subset of $\mathds{R}^d$, $d\geq 1$. Among the advantages of this
approach is the possibility to study the system at both local and
global levels, where the local structure is determined by the
interaction between the entities dependent on their spatial
positions. A paramount task of this study is to describe how does
the local structure of a given system affects its global behavior.
Typically, the entities interacting with a given entity lie in a
compact subset of $\mathcal{H}$. If $\mathcal{H}$ itself is compact,
the qualitative difference between the global and the local is
inessential. Hence, to see the difference between them one has to
place the system into a noncompact habitat. A finite birth-and-death
system in a noncompact habitat always occupies its compact subset
and  has a tendency to disperse to the empty parts by placing there
the newborn entities. It can thus be classified as a {\it
developing} system in which the possibly existing interactions have
little influence on the global behavior, see \cite{Koz}. Therefore,
the importance of the local competition in determining the global
behavior of a birth-and-death system can fully be understood if the
system is {\it developed}, i.e., is infinite and placed in a
noncompact habitat.

In this work, we continue, cf. \cite{DimaN2,FKKK,KK,Koz}, dealing
with the model introduced in \cite{BP1,Mu}, often called
Bolker-Pacala or Bolker-Pacala-Dieckmann-Law model. Here the habitat
is the Euclidean space $\mathds{R}^d$.  The phase space is the set
$\Gamma$ of all locally finite subsets $\gamma \subset
\mathds{R}^d$, i.e., such that $\gamma_\Lambda:=\gamma\cap\Lambda$
is finite whenever $\Lambda \subset \mathds{R}^d$ is compact. For a
compact $\Lambda$, one defines the counting map $\Gamma \ni \gamma
\mapsto |\gamma_\Lambda|$ ($|\cdot|$ denotes cardinality). Then
$\Gamma$ is equipped with the $\sigma$-field $\mathcal{B}(\Gamma)$
generated by all $\Gamma^{\Lambda,n}:=\{ \gamma \in \Gamma :
|\gamma_\Lambda| = n\}$, $n\in \mathds{N}_0$ and $\Lambda$ compact.
This allows for considering probability measures on $\Gamma$ as
states of the system, including Poisson states in which the entities
are independently distributed over $\mathds{R}^d$.  For the
homogeneous Poisson state $\pi_\varkappa$ with density $\varkappa
>0$ and every compact $\Lambda$, one has
\begin{equation}
  \label{J1}
\pi_\varkappa (\Gamma^{\Lambda,n}) = \left(\varkappa {\rm
V}(\Lambda)\right)^n \exp\left( - \varkappa {\rm V}(\Lambda)
\right)/ n!, \qquad n \in \mathds{N}_0,
\end{equation}
where ${\rm V}(\Lambda)$ denotes Lebesgue's measure (volume) of
$\Lambda$. A state $\mu$ can be called \emph{sub-Poissonian}  if for
each compact $\Lambda \subset \mathds{R}^d$, the following holds
\begin{equation}
  \label{J2}
\forall n\in \mathds{N}_0 \qquad \mu(\Gamma^{\Lambda,n}) \leq
C_\Lambda \varkappa_\Lambda^n/n!,
\end{equation}
with some positive constants $C_\Lambda$ and $\varkappa_\Lambda$.
Thus, the sub-Poissonian states are characterized by the lack of
\emph{heavy tails} or \emph{clustering}. The entities in such a
state are either independent in taking their positions or `prefer'
to stay away of each other. The set of \emph{finite} configurations
$\Gamma_0:= \cup_{n\in \mathds{N}_0}\{ \gamma \in \Gamma:
|\gamma|=n\}$ is clearly measurable. In a state with the property
$\mu(\Gamma_0)=1$, the system  is ($\mu$-almost surely)
\emph{finite}.   In this note, we consider infinite systems and
hence deal with states $\mu$ such that $\mu(\Gamma_0)=0$.

In dealing with states on $\Gamma$ one employs {\it observables} --
appropriate functions $F:\Gamma \rightarrow \mathds{R}$. Their
evolution is obtained from the Kolmogorov equation
\begin{equation}
 \label{R2}
\frac{d}{dt} F_t = L F_t , \qquad F_t|_{t=0} = F_0, \qquad t>0,
\end{equation}
where the generator $L$ specifies the model. The states' evolution
is then obtained from the Fokker--Planck equation
\begin{equation}
 \label{R1}
\frac{d}{dt} \mu_t = L^* \mu_t, \qquad \mu_t|_{t=0} = \mu_0,
\end{equation}
related to that in (\ref{R2}) by the duality $\mu_t(F_0) =
\mu_0(F_t):= \int_{\Gamma} F_t(\gamma) \mu_0(d \gamma)$. The model
discussed in this work is specified by the following
\begin{eqnarray}
 \label{L}
 \left(L F \right)(\gamma) & = & \sum_{x\in \gamma} E^{-} (x, \gamma \setminus x) \left[F(\gamma \setminus x) - F(\gamma) \right]\\[.2cm]
 & + & \int_{\mathds{R}^d} E^{+} (x, \gamma ) \left[F(\gamma \cup x) - F(\gamma) \right]dx, \nonumber
\end{eqnarray}
where $E^{+} (x, \gamma )$ and $E^{-} (x, \gamma )$ are
state-dependent birth  and death rates, respectively. We take them
in the following forms
\begin{equation}
 \label{J4}
E^{+} (x, \gamma ) = \sum_{y\in \gamma} a^{+} (x-y),
\end{equation}
\begin{equation}
 \label{J6}
 E^{-} (x, \gamma ) =  m + \sum_{y\in \gamma} a^{-} (x-y),
\end{equation}
where $a^{+}\geq 0$ and $a^{-}\geq 0$ are the \emph{dispersal} and
\emph{competition kernels}, respectively, $m\geq 0$  is the
intrinsic mortality rate.  This model plays a significant role in
the mathematical theory of ecological systems, see \cite{O}. Its
recent study can be found in \cite{DimaN2,FKKK,KK}. The particular
case of (\ref{L}), (\ref{J6}) with $a^{-} \equiv 0$ is the {\it
continuum contact model}, see  \cite{KKP} and the references
therein. In this work,
 we aim at understanding the ecological
consequences of the competition presented in (\ref{L}).
\begin{remark}
  \label{Alt}
\vskip.1cm For the kernels $a^{\pm}$, one has the following
possibilities:
\begin{itemize}
  \item[(i)] (\emph{short dispersal}) there exists $\theta>0$ such that $a^{-} (x) \geq \theta
  a^{+}(x)$ for all $x\in
  \mathds {R}^d$;
\item[(ii)] (\emph{long dispersal}) for each $\theta>0$, there exists  $x\in
  \mathds {R}^d$ such that $a^{-} (x) < \theta
  a^{+}(x)$.
\end{itemize}
\end{remark}
In case (i), $a^{+}$ decays faster than $a^{-}$, and hence each
daughter entity competes with  her mother. Such models are usually
employed to describe the dynamics of cell communities, see
\cite{DimaRR}. An instance of the short dispersal is given by
$a^{+}$ with finite range, i.e., $a^{+}(x)\equiv 0$ for all $|x|
\geq r$, and $a^{-}(x)>0$ for such $x$. In case (ii), $a^{-}$ decays
faster than $a^{+}$, and hence some of the offsprings can be out of
reach of their parents. Models of this kind can be adequate, e.g.,
in plant ecology with the long-range dispersal of seeds. In this
study, the model parameters are supposed to satisfy the following.
\begin{assumption}
  \label{Ass1}
The kernels $a^{\pm}$ in (\ref{J4}) and (\ref{J6}) are continuous
and belong to $L^1 (\mathds{R}^d) \cap L^\infty (\mathds{R}^d)$.
According to this we set $\langle a^{\pm}\rangle =
\int_{\mathds{R}^d} a^{\pm } (x ) dx$ and $\|a^{\pm}\| = \sup_{x\in
\mathds{R}^d} a^{\pm}(x)$.
\end{assumption}
Like in \cite{FKKK,KK}, the evolution of states will be described by
means of correlation functions. To explain the essence of this
approach let us consider the set of all compactly supported
continuous functions $\theta:\mathbb{R}^d\to (-1,0]$. For a state,
$\mu$, its {\it Bogoliubov} functional is
\begin{equation}
  \label{I1}
B_\mu (\theta) = \int_{\Gamma} \prod_{x\in \gamma} ( 1 + \theta (x))
\mu( d \gamma),
\end{equation}
with $\theta$ running through the mentioned set of functions. For
the homogeneous Poisson measure $\pi_\varkappa$,  it takes the form
\begin{equation*}
B_{\pi_\varkappa} (\theta) = \exp\left(\varkappa
\int_{\mathbb{R}^d}\theta (x) d x \right).
\end{equation*}
Having this in mind we will consider those states $\mu$ for which
the functional (\ref{I1}) can be written down in the form
\begin{eqnarray}
  \label{I3}
B_\mu(\theta) = 1+ \sum_{n=1}^\infty
\frac{1}{n!}\int_{(\mathbb{R}^d)^n} k_\mu^{(n)} (x_1 , \dots , x_n)
\theta (x_1) \cdots \theta (x_n) d x_1 \cdots d x_n,
\end{eqnarray}
where $k_\mu^{(n)}$ is a symmetric element of $L^\infty
((\mathbb{R}^d)^n)$. It is the \emph{$n$-th order correlation
function of} $\mu$. In the contact model,  for each $t>0$ and $n\in
\mathbb{N}$ the correlation functions satisfy the following
estimates
\begin{equation}
  \label{IJ33}
{\rm const}\cdot n! c_t^n \leq k^{(n)}_t(x_1, \dots , x_n) \leq {\rm
const}\cdot n! C_t^n.
\end{equation}
Thus, the corresponding state does not satisfy (\ref{J2}), and hence
the clustering does occur in this model. In view of this, the main
question arising here is whether the competition contained in $L$
can suppress clustering. In such a case, one can say that the
\emph{local self-regulation} takes place in this model. The answer
given Theorems \ref{1tm} and \ref{2tm}  below is in affirmative.

\section{The  Results}

By $\mathcal{B}(\mathds{R})$  and $\mathcal{P}(\Gamma)$ we denote
the sets of all Borel subsets of $\mathds{R}$ and   the set of all
probability measures on $(\Gamma, \mathcal{B}(\Gamma))$,
respectively. Ny definition, the subset $\mathcal{P}_{\rm
exp}(\Gamma)\subset \mathcal{P}(\Gamma)$ consists of all those $\mu$
for which $B_\mu$ can be continued, as a function of $\theta$, to an
exponential type entire function on $L^1 (\mathbb{R}^d)$. It can be
shown that a given $\mu$ belongs to $\mathcal{P}_{\rm exp}(\Gamma)$
if and only if $B_\mu$ can be written down as in (\ref{I3}) where
the correlation functions satisfy
\begin{equation}
\label{I4}
  \|k^{(n)}_\mu \|_{L^\infty
((\mathbb{R}^d)^n)} \leq C \exp( \vartheta n), \qquad n\in
\mathbb{N}_0,
\end{equation}
with some $C>0$ and $\vartheta \in \mathbb{R}$. In other words,
$k^{(n)}_\mu$ satisfies the {\em Ruelle bound}, see \cite[Section
6]{Tobi}. In view of (\ref{I4}), each $\mu\in \mathcal{P}_{\rm
exp}(\Gamma)$ satisfies (\ref{J2}) and hence is sub-Poissonian.

 A function $G:\Gamma_0 \subset \Gamma \to \mathds{R}$ is
$\mathcal{B}(\Gamma)/\mathcal{B}(\mathds{R} )$-measurable, see
\cite{FKKK}, if and only if, for each $n\in \mathds{N}$, there
exists a symmetric Borel function $G^{(n)}: (\mathds{R}^{d})^{n} \to
\mathds{R}$ such that
\begin{equation}
 \label{7}
 G(\eta) = G^{(n)} ( x_1, \dots , x_{n}), \quad {\rm for} \ \eta = \{ x_1, \dots , x_{n}\}.
\end{equation}
 Like in (\ref{7}), we introduce $k_\mu : \Gamma_0
\to \mathds{R}$ such that $k_\mu(\eta) = k^{(n)}_\mu (x_1, \dots ,
x_n)$ for $\eta = \{x_1, \dots , x_n\}$, $n\in \mathds{N}$. We also
set $k_\mu(\emptyset)=1$. Then we pass from (\ref{R1}) to the
corresponding Cauchy problem for the correlation functions
\begin{equation}
  \label{J7}
\frac{d}{dt} k_t = L^\Delta k_t, \qquad k_t|_{t=0} = k_{\mu_0}.
\end{equation}
In non-equilibrium statistical mechanics, the corresponding problem
is known as the {\em BBGKY hierarchy}. The action of $L^\Delta$
presents as follows, cf. \cite{DimaN2,FKKK,KK},
\begin{eqnarray}
  \label{J8}
\left( L^\Delta k \right) (\eta)& = & (L^{\Delta,-}k)(\eta) + \sum_{x\in\eta} E^{+} (x , \eta\setminus x) k(\eta \setminus x)\\[.2cm] \nonumber
& + & \int_{\mathds{R}^d} \sum_{x\in \eta} a^{+} (x-y) k(\eta
\setminus x \cup y) d y, \nonumber
\end{eqnarray}
where
\begin{equation}
 \label{J8a}
(L^{\Delta,-}k)(\eta)  := - E^{-}(\eta) k(\eta) -
\int_{\mathds{R}^d}\left(\sum_{y\in \eta} a^{-}(x-y) \right)
k(\eta\cup x)d x,
 \end{equation}
and
\begin{equation}
  \label{J9}
  E^{-}(\eta) := \sum_{x\in \eta} E^{-}(x ,\eta \setminus x) = m|\eta| + \sum_{x\in \eta}\sum_{y\in \eta\setminus x}a^{-} (x-y) .
\end{equation}
By (\ref{I4}) it follows that $\mu \in \mathcal{P}_{\rm
exp}(\Gamma)$ implies $|k_\mu (\eta)| \leq C \exp( \vartheta
|\eta|)$, holding for $\lambda$-almost all $\eta\in \Gamma_0$, some
$C>0$, and $\vartheta\in \mathds{R}$. In view of this, we set
\begin{equation}
 \label{18}
\mathcal{K}_\vartheta := \{ k:\Gamma_0\to \mathds{R}:
\|k\|_\vartheta <\infty\},
\end{equation}
where $\|k\|_\vartheta = \esssup_{\eta \in \Gamma_0}\left\{ |k_\mu
(\eta)| \exp\big{(} - \vartheta
  |\eta| \big{)} \right\}$.
Clearly, (\ref{18})  defines a Banach space. In the following, we
use the ascending scale of such spaces $\mathcal{K}_\vartheta$,
$\vartheta \in \mathds{R}$, with the property $\mathcal{K}_\vartheta
\hookrightarrow \mathcal{K}_{\vartheta'}$ for $\vartheta <
\vartheta'$. Here $\hookrightarrow$  denotes continuous embedding.
Then $\mathcal{K} := \cup_{\vartheta \in \mathds{R}}
\mathcal{K}_\vartheta$ is equipped with the corresponding inductive
topology that turns it into a locally convex space.

For each $\vartheta \in \mathds{R}$ and $\vartheta'>\vartheta$, the
expressions in (\ref{J8}), (\ref{J8a}) and (\ref{J9}) can be used to
define the corresponding bounded linear operators
$L^\Delta_{\vartheta' \vartheta}$ acting from
$\mathcal{K}_\vartheta$ to $\mathcal{K}_{\vartheta'}$. Their
operator norms can be estimated similarly as in \cite[eqs. (3.11),
(3.13)]{KK}, which yields
\begin{equation}
  \label{J11}
\|L^\Delta_{\vartheta' \vartheta}\| \leq \frac{4 (\| a^{+} \| +\|
a^{-} \|  )}{e^2 (\vartheta' - \vartheta)^2} + \frac{ \langle a^{+}
\rangle +m  + \langle a^{-} \rangle
 e^{\vartheta'}}{e (\vartheta' - \vartheta)}.
 \end{equation}
By means of the collection $\{L^\Delta_{\vartheta' \vartheta}\}$
with all $\vartheta\in \mathds{R}$ and $\vartheta'>\vartheta$ we
introduce the corresponding continuous linear operators acting on
$\mathcal{K}$, and thus define the Cauchy problems (\ref{J7}) in
this space. By the (global in time) solutions we mean continuously
differentiable functions $[0,+\infty) \ni t \mapsto k_t \in
\mathcal{K}$ such that both equalities in (\ref{J7}) hold.
\begin{theorem}
  \label{1tm}
Under Assumption \ref{Ass1} the following holds: For each $\mu_0 \in
\mathcal{P}_{\rm exp} (\Gamma)$, the problem in (\ref{J7}) with
$L^\Delta:\mathcal{K}\to \mathcal{K}$ as in (\ref{J8}) -- (\ref{J9})
and (\ref{J11}) has a unique solution $k_t$ such that, for each
$t>0$, there exists a unique state $\mu_t\in \mathcal{P}_{\rm
exp}(\Gamma)$ for which $k_t = k_{\mu_t}$.
\end{theorem}
\begin{theorem}
  \label{2tm}
Let $\vartheta_0$ be such that $k_{\mu_0}\in
\mathcal{K}_{\vartheta_0}$. Then, for all $t\geq 0$, the mentioned
above solution $k_t$,  corresponding to this $k_{\mu_0}$, satisfies
the following estimates:
\begin{itemize}
  \item[{\it (i)}]  {\rm Case $\langle a^{+} \rangle>0$ and $m\in [0,\langle a^{+}
\rangle]$:} for each $\delta <m$ (long dispersal) or $\delta \leq m$
(short dispersal), there exists a positive $C_\delta$ such that
$\log C_\delta \geq \vartheta_0$ and
\begin{equation*}
k_t (\eta) \leq C_\delta^{|\eta|} \exp  \bigg{(}  (\langle a^{+}
\rangle - \delta )|\eta| t \bigg{)}, \qquad \eta \in \Gamma_0.
\end{equation*}
  \item[{\it (ii)}] {\rm Case $\langle a^{+} \rangle>0$ and $m>\langle a^{+}
\rangle$:} for each $\varepsilon \in (0,m-\langle a^{+} \rangle)$,
there exists a positive $C_\varepsilon$ such that $\log
C_\varepsilon \geq \vartheta_0$ and
\begin{equation}
  \label{Bedd}
k_t (\eta) \leq C_\varepsilon^{|\eta|} \exp (-\varepsilon t) ,
\qquad \eta \neq \emptyset.
\end{equation}
\item[{\it (iii)}]{\rm Case $\langle a^{+} \rangle =0$:}
\begin{equation}
  \label{Est}
  k_t (\eta) \leq  k_0(\eta) \exp \left[ - E^{-}(\eta) t \right], \qquad \eta \in \Gamma_0.
\end{equation}
\end{itemize}
If $m=0$ and $a^{-} (x) = \theta a^{+} (x)$, then
\begin{equation}
  \label{bed1}
k_t (\eta) = \theta^{-|\eta|}, \qquad t\geq 0,
\end{equation}
is a stationary solution.
\end{theorem}
The proof of these statements is based on the following
\begin{lemma}
  \label{J1lm}
Let $a^{\pm}$ satisfy Assumption \ref{Ass1}. Then one finds $b \geq
0$ and $\theta >0$ such that
\begin{equation}
  \label{J13}
b |\eta| + \sum_{x\in \eta}\sum_{y\in \eta\setminus x} a^{-}(x-y)
\geq \theta \sum_{x\in \eta}\sum_{y\in \eta\setminus x} a^{+}(x-y),
\end{equation}
holding for all $\eta\in \Gamma_0$.
\end{lemma}
The proof of the lemma (quite technical) can be found in
\cite{KKPr}.

\section{Comments and comparison}

The condition of continuity of $a^{-}$ in Assumption \ref{Ass1} can
be relaxed. In fact,  it is enough to assume that $a^{-}$ is
measurable and separated away from zero in some ball. For $a^{+}$,
it is enough to have a continuous $\tilde{a}^{+}\in L^\infty
(\mathds{R}^d)\cap L^1(\mathds{R}^d)$ such that $\tilde{a}^{+}(x)
\geq a^{+}(x)$ for almost all $x$.

By Theorem \ref{1tm},  adding competition to the continuum contact
model, cf. (\ref{IJ33}),  yields the local self-regulation -- no
matter how long the dispersal is. In the short dispersal case,  the
inequality in (\ref{J13}) readily holds with $b =0$. Then the most
intriguing question here is whether it can hold in the long
dispersal case. In \cite[Proposition 3.7]{KK}, it was shown that
measurable $a^{+}$ and $a^{-}$ satisfy (\ref{J13}) with some $b$ and
$\theta$ if $a^{-}(x)$ is separated away from zero for $|x|<r$ with
some $r>0$, and $a^{+}(x)\equiv 0$ for $|x| \geq R$ with some $R>0$
with the possibility $R>r$. Another choice of $a^{+}$ and $a^{-}$
satisfying (\ref{J13}) can be, see \cite[Proposition 3.8]{KK},
\[
a^{\pm} (x) = \frac{c_{\pm}}{( 2 \pi\sigma_{\pm}^2)^d/2} \exp\left(
- \frac{1}{2 \sigma_{\pm}^2} |x|^2\right),
\]
with all possible  $c_{\pm} >0$ and $\sigma_{\pm} >0$. An important
example of $a^{\pm}$ which both Propositions 3.7 and 3.8 of
\cite{KK} do not cover is the case where $a^{-}$ has finite range
and $a^{+}$ is Gaussian as above. The novelty of our present --
rather unexpected -- result is that (\ref{J13}) is satisfied
\emph{for any} $a^{+}$ and $a^{-}$ as in Assumption \ref{Ass1}, and
hence the local self-regulation is achieved by applying \emph{any
kind of competition.}

Now let us compare our results  with those of \cite{DimaN2,FKKK}. In
\cite{DimaN2}, the model was supposed to satisfy the conditions, see
\cite[Eqs. (3.38) and (3.39)]{DimaN2}, which  can be formulated as
follows: (a) condition (i) in Remark \ref{Alt} holds with a given
$\theta>0$; (b) $m
> 16 \langle a^{-} \rangle / \theta$  for
this $\theta$. Then the global evolution $k_0 \mapsto k_t$  was
obtained in $\mathcal{K}_\vartheta$ with some $\vartheta\in
\mathds{R}$ by means of a $C_0$-semigroup. No information was
available on whether $k_t$ is a correlation function and hence on
the sign of $k_t$. In \cite{FKKK}, the restrictions were relaxed to
imposing the short dispersal condition. Then the evolution $k_0
\mapsto k_t$ was obtained in a scale of Banach spaces
$\mathcal{K}_\alpha$ as in Theorem \ref{1tm}, but on a bounded time
interval. Like  in \cite{DimaN2}, also here no information was
obtained on whether $k_t$ is a correlation function.

Theorem \ref{2tm} gives a complete characterization of the evolution
$k_0 \mapsto k_t$. For $m < \langle a^{+} \rangle$ (short dispersal)
or $m \leq \langle a^{+} \rangle$ (long dispersal), the evolution
described in Theorem \ref{1tm} takes place in an ascending scale
$\{\mathcal{K}_{\vartheta_t}\}_{t\geq 0}$ of Banach spaces. If $m
> \langle a^{+} \rangle$,  the
evolution holds in one and the same space. The only difference
between the cases of long and short dispersals is that one can take
$\delta =m$ in the latter case. This yields different results for $m
= \langle a^{+} \rangle$, where the evolution takes place in the
same space $\mathcal{K}_\vartheta$ with $\vartheta = \log C_m$. Note
also that for $m=0$, one should take $\delta<0$. For $m> \langle
a^{+} \rangle$, it follows from (\ref{Bedd}) that the population
dies out: for $\langle a^{+} \rangle >0$, the following holds
\[
k^{(n)}_{\mu_t} (x_1 ,\dots , x_n) \leq e^{-\varepsilon
t}k^{(n)}_{\mu_0} (x_1 ,\dots , x_n), \quad t>0,
\]
for some $\varepsilon \in (0, m - \langle a^{+} \rangle)$, almost
all $(x_1 ,\dots , x_n)$, and each $n\in \mathbb{N}$.
 For $m>0$ and $\langle a^{+} \rangle =0$, by
(\ref{Est}) we get
\begin{equation*}
 k^{(n)}_{\mu_t} (x_1 ,\dots , x_n) \leq \exp\left(-  n m t \right) k^{(n)}_{\mu_0} (x_1 ,\dots , x_n), \quad t>0.
\end{equation*}
This means that  $ k^{(n)}_{\mu_t} (x_1 ,\dots , x_n)\to 0$ as $n\to
+\infty$ for sufficiently big $t>0$. This phenomenon does not follow
from (\ref{Bedd}). Finally, we mention that (\ref{bed1}) corresponds
to a special case of short dispersal. Until the present work no
results on the extinction as in (\ref{Bedd}) and on the case of
$a^{+} \equiv 0$ were known.
\section*{Acknowledgment}
The present research was supported by the European Commission under
the project STREVCOMS PIRSES-2013-612669 and by the SFB 701
``Spektrale Strukturen and Topologische Methoden in der Mathematik".

\end{document}